\documentclass{article}

\usepackage{amssymb, amsmath, amsthm, bbm}
\input xy
\xyoption{all}
\newtheorem{Theorem}{Theorem}
\newtheorem{Lemma}{Lemma}
\newcommand{\R}{{\mathbb{R}}}

\numberwithin{equation}{section}
\numberwithin{Theorem}{section}
\numberwithin{Lemma}{section}

\begin{document}

\title{An Efficient Approximation of the Traveling Salesman Polytope Using Lifting Methods}
\author{Ellen Veomett}
\date{October, 2006}

\maketitle

\begin{abstract}
For the Traveling Salesman Polytope on $n$ cities $T_n$, we construct its approximation $Q_k, k=1, 2, \dotsc,
n^{1/3}$ using a projection of a polytope
whose number of facets is polynomial in $n$ (of degree
linear in $k$).  We show that
$T_n$ is contained in
$Q_k$ for each $k$, and that the scaling of $Q_k$ by $\frac{k}{n}+O(\frac{1}{n})$ is
contained in $T_n$ for each $k$.  We show that certain facets of $T_n$ lie on the boundary of $Q_k$.

\end{abstract}

\section{Introduction and Results}

For many interesting convex bodies $X$ in a vector space $V$, given a point $x \in V$, the question ``is $x$ in $X$?'' is
 difficult to answer.  This fact has generated work in the direction of finding another set $Y$ which is ``close'' to
$X$ in some way for which the membership question is ``easy'' to answer.  The following results have been the motivation for this paper.

\subsection{Polytope Projection and Successive Approximation}

Given a convex body $X$ in a vector space $V$, a natural type of set to use to approximate $X$ is a polytope. 
Unfortunately, in order to get a ``good'' approximation, this may require that the polytope have exponentially many facets (exponential in $\dim(X)$).  For example, if $B_d$ is the Euclidian ball in $\R^d$, any polytope containing $B_d$ must have exponentially many facets to have its volume be within a factor $c^d$ of the volume of $B_d$ for any contstant $c$ (see, for example, section 13.2 of \cite{Matousek}).

This has led to the idea of approximating convex bodies by \textit{projections} of polytopes, the point being  that the
projection of a polytope may have many more facets than the original polytope.  Ben-Tal and Nemirovski  have exploited
this fact in the case of a Euclidean Ball.  In \cite{BT-Nem}, they proved that for any $\epsilon>0$,  and nonnegative
integer $n$, there exists $N=O(n \log(\epsilon^{-1}))$, a polytope  $P$ with no more than $N$ facets, and a
linear transformation $T$ such that
\begin{equation*}
T(P) \subset B_n \subset (1+\epsilon)T(P)
\end{equation*}
where $B_n = \{x \in \R^n : ||x|| \leq 1\}$ is the Euclidean unit ball in $\R^n$.  Note that an arbitrary convex $X$ can
be  approximated by an ellipsoid within a factor of $\dim(X)$ (see, for example, section V.2 of \cite{B:book}).   Thus,
for $n = \dim(X)$, this automatically gives a polytope $P$ whose number of facets is polynomial in $n$ and
$\log(\epsilon^{-1})$ and a linear transformation $T$ such that
\begin{equation*}
T(P) \subset X \subset (1+\epsilon)nT(P)
\end{equation*}

Sherali and Adams \cite{Sherali-Adams}, Lov\'{a}sz and Schrijver \cite{Lovasz-Schrijver}, and Lasserre  
\cite{Lasserre:2}  have also used projections of polytopes.  In each of these instances, the authors constructed 
successive relaxations of a 0-1 polytope (each of which was a projection of another polytope), such that in the $n$th step,
 the 0-1 polytope is achieved: $P = K^n \subset K^{n-1} \subset \dotsb \subset K^1 \subset K$.  Metric properties of
these sets are not known.  For specifics, 
as well as a comparison of the methods, see \cite{Laurent}.

\subsection{Approximation of the STSP}\label{construction}

The Symmetric Traveling Salesman Polytope (STSP) can be described as follows: recall that a Hamiltonian cycle in the complete graph on $n$ vertices $K_n$ is a cycle which visits every vertex exactly once.  To each Hamiltonian cycle in $K_n$, we can associate its incidence matrix $A = (a_{ij})$ where
\begin{equation*}
a_{ij}= \begin{cases}
1 & \text{if the cycle contains edge } \{i,j\} \\
0 & \text{if the cycle does not contain edge } \{i,j\}
\end{cases}
\end{equation*}
It is called symmetric because there is a similar notion in the case of a digraph (a graph where the edges have an orientation),
 which is the Asymmetric Traveling Salesman polytope.

Note that each matrix corresponding to a Hamiltonian cycle is a symmetric 0-1 matrix in $\R^{n^2}$ with 0s on the
diagonal.  Given a particular matrix corresponding to a Hamiltonian cycle, any other such matrix can be achieved from it
by  simultaneously permuting rows and columns (this corresponds to permuting the labels on the vertices of the graph).
  The Traveling Salesman Polytope is the convex hull of all adjacency matrices corresponding to Hamiltonian cycles in $K_n$.
  Note that the vertices of the STSP are matrices which correspond to cycles.  To each cycle we can associate  a
permutation of the numbers $\{1, 2, \dotsc, n\}$ beginning with the number 1 where the permutations  $(1, m_2, m_3,
\dotsc, m_n)$ and $(1, m_n, m_{n-1}, \dotsc, m_2)$ are identified.  We will use the descriptions  of the vertices as
matrices, cycles, and permutations interchangeably.  With the permutation description of a Hamiltonian cycle,  it is
not hard to see that there are $\frac{(n-1)!}{2}$ different Hamiltonian cycles in $K_n$.

The STSP has been studied widely, though a complete description of it via linear
inequalities is unknown.  It
is clearly not full dimensional in $\R^{n^2}$, being the convex hull of symmetric matrices with 0s on the diagonal.  
It is not hard to show that its dimension is $\frac{n(n-3)}{2}$.  For more information on the STSP and the associated 
Traveling Salesman Problem, see, for example, Chapter 58 of \cite{Schrijver:book}.   Linear optimization over the STSP
and the membership question for the STSP are known to be NP-hard.

Let $X$ denote the set of matrices corresponding to Hamiltonian cycles in $K_n$.  Instead of working directly with the STSP, we will be working with its polar.  Note that the barycenter of the STSP is the matrix $Z = (z_{ij})$ where
\begin{equation*}
z_{ij} = \begin{cases} \frac{2}{n-1} & \text{if } i \not= j \\
0 & \text{if } i=j
\end{cases}
\end{equation*}
Since the STSP is not full dimensional, in order to get a bounded polar we move $Z$ to the origin, which forces the average value to be 0:
\begin{equation*}
A = \biggl\{\text{linear functions } f: \quad \frac{1}{|X|}\sum_{x \in X} f(x) = 0,\quad  f(x) \leq 1 \quad \text{ for all } x \in X\biggr\}
\end{equation*}
Then for pure convenience, we will reflect $A \mapsto -A$ and then shift $f \mapsto f+1$ to obtain the following description of the dual that will be the one with which we work:
\begin{equation*}
 \begin{split}Q = \biggl\{ \text{functions } f : & f \text{ is linear, that is } f(x) = \left<c,x\right> \text{ for some matrix } c \\
& \text{non-negative, that is } f(x) \geq 0 \text{ for all } x \in X \\
& \text{has average 1, that is } \frac{1}{|X|} \sum_{x \in X} f(x) = 1 \biggr\}
\end{split} 
\end{equation*}
Note that the center of $Q$ is the all ones function: $\mathbbm{1}(x) = 1$ for all $x \in X$.

Since we want to approximate the STSP with a projection of a polytope having not too many facets, this  equates to
approximating the set $Q$ above with a section of a polytope having not too many vertices.   We will view $Q$ as living
in the space $L$ of all linear functions $f:X \to \R$, that is, restrictions of linear functions on $\R^{n^2}$ to $X$.
We in turn view $L$ as living in $\R^X$, the space of
\textit{all} functions $f:X \to
\R$.

With this setup, we have the following:
\begin{Theorem}\label{project1}
Let $n \geq 4$ be an integer and $k \leq n^{1/3}$ be an integer.  There exists a polytope $P_k\subset \R^X$ with $O(n^{4k})$ vertices and a constant 
\begin{equation*}
c_k = \frac{k}{n} + O\bigg(\frac{1}{n}\bigg)
\end{equation*}
such that 
\begin{equation*}
c_k(Q-\mathbbm{1}) \subset P_k\cap L -\mathbbm{1} \subset Q-\mathbbm{1}
\end{equation*}
\end{Theorem}

Note that the dimension of the convex hull of $m$ vertices is $\leq m-1$.  Recall that we denoted by $Z$ the center of the STSP.  Thus, from the remarks above, this immediately gives as a corollary the following:
\begin{Theorem}\label{project2}
Let $T_n$ be the Symmetric Traveling Salesman Polytope and let $k \leq n^{1/3}$ be an integer.  Then there exists 
a  polytope $P_k^\circ$ with $O(n^{4k})$ facets, a linear transformation $T$, and a constant
\begin{equation*}
c_k = \frac{k}{n} + O\bigg(\frac{1}{n}\bigg)
\end{equation*}
such that for $Q_k = T(P_k^\circ)$ we have
\begin{equation*}
c_k(Q_k-Z) \subset T_n -Z\subset Q_k-Z
\end{equation*}
\end{Theorem}

One notable aspect of this approximation is that the scaling factor gives us a \textit{metric} bound on how far our
approximating set can be from the STSP.

\subsection{Computability Remarks} 
Say that $P$ is a polytope living in $\R^N$ with $O(N)$ facets and $\pi: \R^N \to \R^n$ is a projection.   Then
deciding if a point is in $\pi(P)$ becomes a linear programming problem in $O(N)$ equations and variables.   Linear
programming is decidable in time polynomial in the number of equations and variables (see, for example Chapter 3 of \cite{GLS}).

\section{Projection Construction}\label{ProjectionConstruction}

Recall that we consider $Q$ to be a subset of $\R^X$, and we will be approximating $Q$ by a polytope in $\R^X$ with not too many vertices.  Thus, we need to describe which functions in $\R^X$ will serve as our vertices (the functions of which we will take the convex hull).  Fixing a $k < n^{1/3}$ (the reason for this restriction will be evident later), we will consider functions, each indexed by a particular  subset of the edges of the complete graph $K_n$.  We only consider subsets of edges of $K_n$ which could correspond to a subset of a Hamiltonian cycle in $K_n$; namely, a subset of edges which correspond to disjoint paths.  We call such subsets ``path subsets.'' 

Given a path subset $\Gamma$ with $k$ edges in it, note that the lengths of the disjoint paths in $\Gamma$ are a partition $\pi$ of $k$.  We call this partition the ``partition type'' of $\Gamma$.  From the following Lemma, we can see that the number of Hamiltonian cycles containing all edges in $\Gamma$ depends only on $k$ and on the number of parts in
$\pi$ (i.e., on the partition type of $\Gamma$).

\begin{Lemma}\label{paths}
Let $(k_1, k_2, \dotsc, k_m)$ be a partition of $k\leq n-1$ ($k+m \leq n$) and $K_n $ the complete graph on $n$ vertices.  Let $p_1, p_2, \dotsc, p_m$ be disjoint paths in $K_n$ of length $k_1, \dotsc, k_m$ respectively.  Then the number of Hamiltonian cycles in $K_n$ containing all of paths $p_1, \dotsc, p_m$ is:
\begin{equation*}
2^{m-1}(n-k-1)!
\end{equation*}
\end{Lemma}

\begin{proof}

Note that the restriction $k+m \leq n$ assures that it is possible to find disjoint paths in $K_n$ of lengths $k_1, \dots, k_m$.  Any cycle containing the paths $p_1, \dotsc, p_m$ can be written uniquely as a sequence of numbers, beginning with path $p_1$ in a particular orientation.  Thinking of the remaining paths as blocks with 2 orientations and the remaining numbers as blocks with a single orientation, we find that there are $2^{m-1}(n-k-1)!$ ways of ordering and orienting the remaining blocks.  Each of these orders and orientations corresponds uniquely to a Hamiltonian cycle containing paths $p_1, \dotsc, p_m$.
\end{proof}

Thus, we shall denote by $a_{\pi(\Gamma)} = a_\pi$  the number of Hamiltonian cycles containing all edges in $\Gamma$, where $\Gamma$ has partition type $\pi$. 

First we define $[A]$ to be the indicator function of a set $A$.  That is,
\begin{equation*}
[A](x) = \begin{cases}
1 & \text{ if } x \in A \\
0 & \text{ if } x \not\in A
\end{cases}
\end{equation*}
By abuse of notation, for a single point $x$ we write $[x]$ for its indicator function instead of $[\{x\}]$.

Now we can define

\begin{equation*}
g_\Gamma = \frac{|X|}{a_{\pi(\Gamma)}} \sum_{\substack{x \in X \\ \Gamma \subset x} } [x]
\end{equation*}
Thus, 
\begin{equation*}
g_\Gamma(x) = \begin{cases} 
\frac{|X|}{a_{\pi(\Gamma)}} & \text{ if }  x  \text{ contains all edges in } \Gamma \\
0  & \text{ otherwise} 
\end{cases}
\end{equation*}
Note that the constant is chosen so that $g_\Gamma$ has an average of 1 on $X$.

Given a $k<n^{1/3}$ and a partition $\pi$ of $k$, we will now describe a linear operator 
\begin{equation*}
T_\pi : \R^X \to \text{span}\{g_\Gamma: \Gamma \text{ is a path subset}, |\Gamma| = k, \Gamma \text{ has partition type } \pi\}
\end{equation*}
Define 
\begin{align*}
A_{\Gamma} &= \{x \in X : x \text{ contains all the edges in } \Gamma\} \\
B_\pi &= \{\Gamma \text{ a path subset}: \Gamma \text{ has partition type } \pi\}
\end{align*}
Note that $|A_{\Gamma}| = a_{\pi(\Gamma)} = a_\pi$.  For a $f \in \R^X$, we define 

\begin{equation*}
T_\pi(f) = 
\frac{1}{|B_\pi|} \sum_{\Gamma \in B_\pi} \alpha_\Gamma g_\Gamma
\end{equation*}

where 

\begin{equation*}
\alpha_\Gamma = \frac{1}{a_{\pi(\Gamma)}} \sum_{x \in A_\Gamma} f(x)
\end{equation*}

In words, the operator $T_\pi$ is  a weighted sum of the functions $g_\Gamma$ where $\Gamma$ has partition type $\pi$.  The weight of a particular $g_\Gamma$ is \begin{equation*}
\frac{\text{average value of } f \text{ on } x \text{ containing } \Gamma}{\text{number of }g_\Gamma \text{ with partition type } \pi}
\end{equation*}

Suppose that $f$ is a function which has average value 1 on $X$.  Then the sum of the above coefficients is:

\begin{equation*}\begin{split}
\frac{1}{|B_\pi|} \sum_{\Gamma \in B_\pi} \frac{1}{a_{\pi(\Gamma)}} \sum_{x \in A_{\pi(\Gamma)}} f(x) &\\
=& \sum_{x \in X} \frac{|\{\text{path subsets in } x \text{ of partition type } \pi\}|}{|B_\pi| a_\pi}  f(x) \\
=& \frac{1}{|X|} \sum_{x \in X} f(x) = 1
\end{split}\end{equation*}

Thus, we can see that if $f$ is a function which has average value 1 on $X$, then $T_\pi(f)$ is a convex combination of the $g_\Gamma$s.  Our goal is to understand how $T_\pi$ acts on the $g_{st}$ defined as follows:

\begin{equation*}
g_{st} = \frac{n-1}{2} \sum_{x \; \text{contains edge} \; \{s,t\}} [x]
\end{equation*}
so that 
\begin{equation*}
g_{st}(x) = \begin{cases}
 \frac{(n-1)}{2}  & \text{ if }  x  \text{ contains edge } \{s,t\} \\
 0  & \text{ otherwise}
\end{cases}
\end{equation*}
Note that $g_{st}$ is the particular case of $g_\Gamma$ where $\Gamma$ consists of a single edge.

It is clear that $Q$ (defined in section \ref{construction}) is contained in the affine span of the $g_{st}$.  Our agenda
at this point is to obtain a convex
combination $T$ of the linear maps $T_\pi$ and a constant $c_k = \frac{k}{n} + O\left(\frac{1}{n}\right)$ such that
\begin{equation}\label{bigT}
T(g_{st}) =(1-c_k)\mathbbm{1} + (c_k)  g_{st}
\end{equation}
The linear map $T$ will  act the same way on the affine span of the $g_{st}$.   We will define $P_k $ to be the
convex hull of of the $g_\Gamma$ used in $T$, and $L$ to be the subspace of linear functions on $X$.  Because each of
the $g_\Gamma$ are nonnegative on $X$ and have average value 1 on $X$, we can see that any function in $P_k \cap L$
will also be in $Q$.  And since, as mentioned already, $Q$ is contained in the affine span of the $g_{st}$ and $T$ will
act the same way on the affine span of the $g_{st}$, we will have:
\begin{equation*}
c_k(Q -\mathbbm{1}) \subset P_k \cap L-\mathbbm{1} \subset Q-\mathbbm{1}
\end{equation*}
which will give us Theorem \ref{project1}.   Thus, we proceed in finding such a linear map $T$.

Let $\pi$ be a partition of $k<n^{1/3}$.  Then
\begin{equation}\label{T(g_{st})}
T_\pi(g_{st}) = \frac{1}{|B_\pi|} \sum_{\Gamma \in B_\pi} \left( \frac{1}{a_{\pi(\Gamma)}} \sum_{x \in A_\Gamma} g_{st}(x) \right) g_\Gamma 
\end{equation}

For a partition $\pi = (k_1, k_2, \dotsc, k_m)$ of $k$,  we want to calculate exactly what the function $T_\pi(g_{st})$ does.  In equation \eqref{T(g_{st})}, we firstly note that $g_{st}$ takes on only values 0 and $\frac{n-1}{2}$, and that $g_\Gamma$ takes on only values 0 and $\frac{|X|}{a_\pi}$.  Thus, we can see that
\begin{align}
T_\pi(g_{st})(y) =& \frac{|X|(n-1)}{a_\pi 2 |B_\pi| a_\pi}&\# & \{\text{pairs } (x,\Gamma) \text{ such that } x \in A_\Gamma, g_{st}(x) =
\frac{n-1}{2}, \notag \\
& & & g_\Gamma(y) = \frac{|X|}{a_\pi}\} \notag \\
=& \frac{|X|(n-1)}{a_\pi 2 |B_\pi| a_\pi}&\# & \{x: \{s,t\} \in x, \Gamma \subset x \text{ for some } \Gamma \subset y
\notag \\
& & &\text{ of partition type } \pi\} \label{operator}
\end{align}

Note that if $\Gamma$ contains edges corresponding to a path of length $\geq 2$ with $s$ and $t$ as its endpoints, there
cannot be a Hamiltonian cycle containing all
edges in $\Gamma$ as well as the edge $\{s,t\}$ (recall that $k \leq n^{1/3}$).  If  $\Gamma$ contains edges
corresponding to a path where $s$ is connected to two vertices, neither of which is $t$, there cannot be a Hamiltonian cycle
 containing all edges in $\Gamma$ as well as the edge $\{s,t\}$ (and similarly for $t$) because any vertex in a
Hamiltonian cycle has exactly 2 vertices adjacent to it.  There are 4 ways in which a Hamiltonian cycle \textit{can}
contain all edges in $\Gamma$ as well as the edge $\{s,t\}$:
\begin{enumerate}
\item $\{s,t\}$ is an edge in $\Gamma$.
\item $s$ and $t$ are each endpoints of different paths in $\Gamma$
\item Exactly one of $s$ or $t$ is an endpoint of a path in $\Gamma$, and the other does not appear in $\Gamma$
\item Neither $s$ nor $t$ appear in $\Gamma$
\end{enumerate}

Suppose that $\Gamma$ has partition type $\pi = (k_1, \dotsc, k_m)$.  If $\{s,t\} \in \Gamma$, then $\Gamma \cup
\{\{s,t\}\}$  again contains $k$ edges and has a partition type with $m$ parts.  If $s$ and $t$ are each endpoints of separate paths
in  $\Gamma$, then  $\Gamma \cup \{\{s,t\}\}$ contains $k+1$ edges and has a partition type with $m-1$ parts.  If exactly
one  of $s$ or $t$ is an endpoint of a path in $\Gamma$, then $\Gamma \cup \{\{s,t\}\}$ contains $k+1$ edges and has a
partition  type with $m$ parts.  And if neither $s$ nor $t$ appear in $\Gamma$, then $\Gamma \cup \{\{s,t\}\}$ contains
$k+1$  edges and has a partition type with $m+1$ parts.

Using Lemma \ref{paths} and equation \eqref{operator}, we can now more explicitly describe the function $T_\pi(g_{st})$.

\begin{align}\label{explicitT}
T_\pi(g_{st})(y)  = &\frac{|X|(n-1)}{a_\pi 2 |B_\pi| a_\pi} \biggl(|\{\Gamma \subset y : \Gamma \text{ has partition } \pi, \text{ and } \{s,t\} \in \Gamma\}| \notag \\
& \cdot 2^{m-1}(n-k-1)! \notag \\
+& |\{\Gamma \subset y : \Gamma \text{ has partition } \pi, \text{ and } s,t \text{ are each endpoints}  \notag \\
& \text{ of different paths in } \Gamma \}| \cdot 2^{m-2}(n-k-2)! \notag \\
+& |\{\Gamma \subset y : \Gamma \text{ has partition } \pi, \text{ and exactly one of } s,t  \text{ is an endpoint  } \notag \\ 
& \text{ of a path in } \Gamma \}| \cdot 2^{m-1}(n-k-2)! \notag \\
+& |\{ \Gamma \subset y : \Gamma \text{ has partition } \pi \text{ and } s,t  \text{ are not in } \Gamma \}|  \notag \\
& \cdot 2^m(n-k-2)! \biggr)
\end{align}

Thus, we can see that $T_\pi(g_{st})(y)$ depends only on how $s$ and $t$ sit in the paths of  subsets $\Gamma \subset y$ of partition type $\pi$.
  In other words, the number of vertices between $s$ and $t$ in the Hamiltonian path $y$ determines $T_\pi(g_{st})(y)$.
  The following lemma shows that $T_\pi(g_{st})$ takes on no more than $k+1$ values on $X$:

\begin{Lemma}\label{same} 
Let $\pi$ be a partition of the positive integer $k \leq n^{1/3}$.  Let $m_y$  be the number of vertices
between $s$ and $t$ in the shorter path between $s$ and $t$ in the Hamiltonian cycle $y$.   Then for any $y$ such that
$m_y \geq k$, the value of $T_\pi(g_{st})(y)$ is the same.
\end{Lemma}

Thus, this Lemma implies that the range of values taken on by $T_\pi(g_{st})$ can be found by evaluating $T_\pi(g_{st})$ on Hamiltonian cycles having $0, 1, \dotsc, k$ vertices between $s$ and $t$.  The proof of this lemma is postponed until section \ref{proofs}.

We also know how $T_\pi(g_{st})$ acts on $y$ which contain the edge $(s,t)$:

\begin{Lemma}\label{y_0}
Let $\pi$ be a partition of $k\leq n^{1/3}$ and suppose that $y$ contains the edge $\{s,t\}$.  Then $T_\pi(g_{st})(y)= \frac{k+2}{2} + O\left(\frac{1}{n^{2/3}}\right)$.
\end{Lemma}

Again we postpone the proof until section \ref{proofs}.

We note that Lemma \ref{y_0} is key to the fact that it is feasible to find a convex  combination of linear maps $T_\pi$,
resulting in a map $T$ acting as in equation \eqref{bigT}.  What we will show  is that $T_\pi(g_{st})$ is ``almost''
$(1-a)\mathbbm{1}+ag_{st}$ for $a \in (0,1)$.  It isn't \textit{exactly} $(1-a)\mathbbm{1}+ag_{st}$ because $T_\pi(g_{st})$ isn't
 the same value on all Hamiltonian cycles which do not contain the edge $\{s,t\}$.  But the different values it takes on for
Hamiltonian  cycles not containing edge $\{s,t\}$ are very close to each other.   Thus, we take a convex combination of
$T_\pi$s for varying partitions $\pi$ is to ``smooth out'' those differences,  resulting in a single map $T$ taking on
only 2 values: a single value for Hamiltonian cycles containing $\{s,t\}$ and a different  value for Hamiltonian cycles
not containing $\{s,t\}$.  This will imply that $T(g_{st}) = (1-a)\mathbbm{1} + ag_{st}$ for some $a \in (0,1)$ (recall, $T_\pi$
maps  functions with average value 1 to functions with average value 1).  The $T_\pi$s used in the map  $T$ will
involve partitions $\pi$ of varying numbers.  Using Lemma \ref{y_0}, we will find that $T(g_{st})$ will be
$\frac{k+1}{2}(1+O(\frac{1}{n^{1/3}}))$  on Hamiltonian cycles containing the edge $\{s,t\}$.  Thus, 
\begin{equation*}
(1-a)+a\frac{n-1}{2} = \frac{k+1}{2}\biggl(1 + O\biggl(\frac{1}{n^{1/3}}\biggr)\biggr)
\end{equation*}
so that
\begin{equation*}
a = \frac{k}{n}+O\biggl(\frac{1}{n}\biggr)
\end{equation*}

For each of the following Lemmas, we let $y_i$ be a Hamiltonian cycle which has $ i $ vertices between vertex $s$ and vertex $t$, and let $y_{i+1}$ be a Hamiltonian cycle which has $i+1$ vertices between vertex $s$ and vertex $t$.  We let $n_{(m)}$ denote $n(n-1)(n-2)\dotso(n-m+1)$.

\begin{Lemma}\label{change_sequence}
Let $1 \leq i \leq k-1$.  If $\pi$ is the partition $(\underbrace{1,1, \dotsc, 1}_{k \text{ ones}})$ then 
\begin{equation*}
T_\pi(g_{st})(y_{i+1})-T_\pi(g_{st})(y_i) = (-1)^{i+1}\frac{(n-1)(n-2k)(n-2k-1)k_{(i+1)}}{4(n-k-1)n(n-k-1)_{(i+2)}}
\end{equation*}
\end{Lemma}

\begin{Lemma}\label{k,1snakes}
If $\pi$ is the partition $(k-1,1)$ for $k \geq 3$ then 
\begin{equation*}
T_\pi(g_{st})(y_{i+1})-T_\pi(g_{st})(y_i) = \begin{cases}
-\frac{(n-1)(n-k-3)}{n(n-k-1)^2} & 1 \leq i \leq k-3 \\
\frac{3(n-1)}{2n(n-k-1)^2} & i = k-2 \\
\frac{(n-1)}{2n(n-k-1)^2} & i = k-1 \end{cases}
\end{equation*}
\end{Lemma}

\begin{Lemma}\label{ksnake}
If $\pi$ is the partition $(k)$ then 
\begin{equation*}
T_\pi(g_{st})(y_{i+1})-T_\pi(g_{st})(y_i) = \begin{cases}
-\frac{(n-1)}{n(n-k-1)} & 1 \leq i \leq k-2 \\
0 & i = k-1 \end{cases}
\end{equation*}
\end{Lemma}

The proofs of the above Lemmas are postponed to section \ref{proofs}.  Using all of these Lemmas, we can prove Theorem
\ref{project1}.
\begin{proof}[Proof of Theorem \ref{project1}]
Recall from earlier comments that we need only find a convex combination $T$ of the linear maps $T_\pi$  such that
\begin{equation*}
T(g_{st})(y) = \frac{k+1}{2}\bigg(1+O\bigg(\frac{1}{n^{1/3}}\bigg)\bigg)
\end{equation*}
on Hamiltonian cycles $y$ containing the edge $\{s,t\}$, and such that the number of functions $g_\Gamma$ used in $T$
is of order $n^{4k}$.

Recall that we assume $k < n^{1/3}$.  The maps $T_\pi$ which we will use will correspond to the partitions 
\begin{align*}
\pi_{\ell} &= (\ell-1, 1) &  1 \leq \ell \leq 2k \\
\pi_{\ell'} &= (\ell) &   2 \leq \ell \leq 2k+1 \\ 
\pi^* &= (\underbrace{1, 1, \dotsc, 1}_{2k \; 1s}) 
\end{align*}
The way that $T$ is obtained is as follows:  We will use the $T_{\pi_\ell}$ and $T_{\pi_{\ell'}}$ to adjust $T_{\pi^*}$.  Recall (from Lemma \ref{change_sequence}) that if $y_{2k}$ is a Hamiltonian cycle with $2k$ vertices between $s$ and $t$ and $y_{2k-1}$ is a Hamiltonian cycle with $2k-1$ vertices between $s$ and $t$, then 
\begin{multline*}
T_{\pi^*}(g_{st})(y_{2k})-T_{\pi^*}(g_{st})(y_{2k-1}) \\
= (-1)^{2k}\frac{(n-1)(n-2(2k))(n-2(2k)-1)(2k)_{(2k)}}{4(n-2k-1)n(n-2k-1)_{(2k+1)}}
\end{multline*}
Since $2k$ is even, we define $a_{2k} = 0$ and find the positive $b_{2k}$ such that 
\begin{multline*}
T_{\pi^*}(g_{st})(y_{2k})+a_{2k}T_{\pi_{2k}}(g_{st})(y_{2k})+b_{2k}T_{\pi_{2k+1'}}(g_{st})(y_{2k}) \\
-(T_{\pi^*}(g_{st})(y_{2k-1}) +a_{2k}T_{\pi_{2k}}(g_{st})(y_{2k-1})+b_{2k}T_{\pi_{2k+1'}}(g_{st})(y_{2k-1})) = 0
\end{multline*}

From Lemmas \ref{change_sequence} and \ref{ksnake}, we can see that this would imply that $b_{2k} \sim \frac{(2k)_{(2k)}}{(n-2k-1)_{(2k-1)}}$, which
is  negligible if $k < n^{1/3}$.  We can also see that, if $y_{i+1}$ is a Hamiltonian  cycle with $i+1$
vertices between $s$ and $t$ and $y_i$ is a Hamiltonian cycle with $i$ vertices between $s$ and $t$ for $i+1<2k$, then 
\begin{multline*}
T_{\pi^*}(g_{st})(y_{i+1})+a_{2k}T_{\pi_{2k}}(g_{st})(y_{i+1})+b_{2k}T_{\pi_{2k+1'}}(g_{st})(y_{i+1}) \\
-(T_{\pi^*}(g_{st})(y_{i}) +a_{2k}T_{\pi_{2k}}(g_{st})(y_{i})+b_{2k}T_{\pi_{2k+1'}}(g_{st})(y_{i})) 
\end{multline*}
is of the same order as 
\begin{equation*}
T_{\pi^*}(g_{st})(y_{i+1})-T_{\pi^*}(g_{st})(y_{i})
\end{equation*}
In the next step, since $2k-1$ is odd, we define $b_{2k-1} = 0$ and find the positive $a_{2k-1}$ such that 
\begin{multline*}
T_{\pi^*}(g_{st})(y_{2k-1})+a_{2k-1}T_{\pi_{2k-1}}(g_{st})(y_{2k-1})+b_{2k}T_{\pi_{2k+1'}}(g_{st})(y_{2k-1}) \\
-(T_{\pi^*}(g_{st})(y_{2k-2}) +a_{2k-1}T_{\pi_{2k-1}}(g_{st})(y_{2k-2})+b_{2k}T_{\pi_{2k+1'}}(g_{st})(y_{2k-2})) = 0
\end{multline*}
where $y_{2k-1}$ is a Hamiltonian cycle with $2k-1$ vertices between $s$ and $t$, and $y_{2k-2}$ is a Hamiltonian
cycle with $2k-2$ vertices between $s$ and $t$.  From Lemmas \ref{change_sequence} and \ref{k,1snakes}, we can see that
this would imply that $a_{2k-1} \sim \frac{(2k-1)_{(2k-1)}}{(n-(2k-1)-1)_{(2k-3)}}$, which is negligible if $k > 3$ and $k <
n^{1/3}$.   We can also see that, if $y_{i+1}$ is a Hamiltonian cycle with $i+1$ vertices between $s$ and $t$ and $y_i$ is a Hamiltonian cycle
 with $i$ vertices between $s$ and $t$ for $i+1<2k-1$, then 
\begin{multline*}
T_{\pi^*}(g_{st})(y_{i+1})+a_{2k-1}T_{\pi_{2k-1}}(g_{st})(y_{i+1})+b_{2k}T_{\pi_{2k+1'}}(g_{st})(y_{i+1}) \\
-(T_{\pi^*}(g_{st})(y_{i}) +a_{2k-1}T_{\pi_{2k-1}}(g_{st})(y_{i})+b_{2k}T_{\pi_{2k+1'}}(g_{st})(y_{i})) 
\end{multline*}
is of the same order as 
\begin{equation*}
T_{\pi^*}(g_{st})(y_{i+1})-T_{\pi^*}(g_{st})(y_{i})
\end{equation*}

We continue this process, next smoothing out the values between Hamiltonian  cycles having $2k-2$ versus $2k-3$
vertices between $s$ and $t$.  Since $T_{\pi^*}(g_{st})(y_{i+1})-T_{\pi^*}(g_{st})(y_i)$  alternates sign, as we
continue ``smoothing out'' $T_{\pi^*}(g_{st})$, we will alternately use $T_{\pi_\ell}$ and $T_{\pi_{\ell'}}$.  All of the coefficients will be very
small,  except perhaps for the coefficient of $T_{\pi_3}$, which could be up to $\frac{1}{2}$ if $k$ is
close to $n^{1/3}$.  Thus, when we divide by the sum of the coefficients (making a convex combination), the coefficient of $T_{\pi^*}$ will be at least $\frac{1}{2}+O(\frac{1}{n^{1/3}})$.  Thus, from Lemma \ref{y_0}, we know that the final $T$ will have value
$\frac{k+1}{2}(1+O(\frac{1}{n^{(1/3)}}))$ on Hamiltonian cycles containing the edge $\{s,t\}$.  

The number of different functions used in $T$ will be the sum of the number  of different functions used in $T_{\pi^*},
T_{\pi_\ell}$, and $T_{\pi_{\ell'}}$.  The number of different functions  used in $T_{\pi^*}$ is equal to the number of
ways of picking $2k$ disjoint edges from the complete graph $K_n$; i.e. the number  of path subsets of partition type
$(1, 1, \dotsc, 1)$.  To pick $2k$ disjoint edges, we can pick $4k$ numbers from  the set $\{1, 2, \dotsc, n\}$ in
order.  The first two we define as being an ``edge'', the second two we define as being  an ``edge'', etc.  Of course,
we get the same set of edges if we picked two numbers which correspond to an edge in  reverse order (i.e., instead of
picking $i$, then $i+1$ and defining them to be an edge, we picked $i+1$  and then $i$ and defined them to
be an edge).  We also get the same set of edges if two pairs of edges switch places in the  ordering (i.e., instead of
picking in order $i, i+1, i+2, i+3$ and defining edges to be $\{i, i+1\}$  and $\{i+2,
i+3\}$,  we had picked in order $i+2, i+3, i, i+1$).  Thus, we can see that the number of different 
functions used in $T_\pi$ is
\begin{equation*}
\frac{n_{(4k)}}{(2k)!2^{2k}} = O(n^{4k})
\end{equation*}

By similar arguments, we can see that the number of functions used in $T_{\pi_\ell}$ and $T_{\pi_{\ell'}}$ are of order smaller than $n^{4k}$.  Thus, the total number of functions used in $T$ is of order $n^{4k}$, and we have finished the proof.

\end{proof}

\section{Facets on the Boundary}

Although there is no known complete description of the Symmetric Traveling Salesman Polytope as a system of linear inequalities, many facets are known (see, for example, chapter 58 of \cite{Schrijver:book}).  Some well-known facet defining inequalities are the following:

\begin{align}
 0  \leq x_{ij} & \leq 1   &\text{ for each } i,j \label{facet1} \\
 \sum_{\substack{j \in U \\ i \in V-U}} x_{ij} & \geq 2  &\text{ for each } U \subset V \text{ with } \emptyset \not= U \not= V \label{facet2} 
 %\\
 %\sum_{\substack{j \in H  \\ i \not\in H }} x_{ij} + \sum_{\ell=1}^t\sum_{\substack{u \in T_\ell \\ v\not\in T_i}} x_{k\ell}  & \geq 3t+1      \begin{split} & \text{for } H, T_i\subset V, t \geq 3 \text{ odd }, \\
%& H \cap T_i \not= \emptyset \quad \text{for } i=1, 2, \dotsc t  \\  
%& T_i \backslash H \not= \emptyset \quad \text{for } i=1, 2, \dotsc t  \\
%& T_i \cap T_j = \emptyset \quad \text{for } 1 \leq i < j \leq t \end{split} \label{facet3}
\end{align}

A natural question to ask regarding the approximation construction of the previous section would be: which (if any) of the above facets lie on the boundary of the approximating set?  In other words, which (if any) of the inequalities defining our approximating set coincide with one of the above inequalities?

Our construction creates a convex set $P_k$ whose intersection with the space of linear functions lies inside of the dual
of the STSP.   Given our definition of the dual $Q$, we are looking for a function $f \in P_k$ which is linear and for
which the set $\{x \in X: f(x) = 0\}$ is precisely the cycles for which equality holds in one of equations
\eqref{facet1}-\eqref{facet2}.

Fix some $k < n^{1/3}$.  Recall that $P_k$ is the convex hull of functions $g_\Gamma$ which take a single
positive
value on  cycles $x$ containing the edges in $\Gamma$ and 0 on cycles not containing the edges in $\Gamma$.   For the
functions used in $P_k$, $\Gamma$ contains edges which correspond to a path of length $\ell \leq 2k+1$,  or a path of
length $\ell \leq 2k-1$ plus a single disjoint edge, or $2k$ disjoint edges.  

Consider some edge $\{i,j\}$.  Then the set $\Gamma = \{i,j\}$ corresponds to a path of length 1.  We cnsider the function 
 \begin{equation*}
 f_{ij} =  g_\Gamma
 \end{equation*}
Note that $f_{ij} \in P_k$ for $k \geq 1$.  Here we can see that $f_{ij}(x) = 0$ precisely when the cycle $x$ does not contain the edge $\{i,j\}$.  Also, it is clear that $f_{ij}$ takes on the same value for each $x$ containing the edge $\{i,j\}$.  Thus, we can see that the $f_{ij}$s are each linear functions which correspond to the facets defined by the left hand sides of equations \eqref{facet1}
 
 Again we consider some edge $\{i,j\}$.  Now let 
 \begin{equation*}
 X_{ij} = \{\Gamma =\{\{i,a\}, \{i,b\}\}: a,b\in V; a \not= b; a,b\not= j\}
 \end{equation*}
 Note that $X_{ij}$ consists of sets corresponding to paths of length 2.  We define
 \begin{equation*}
 f_{ij}'=\frac{1}{|X_{ij}|}\sum_{\Gamma \in X_{ij}} g_\Gamma
 \end{equation*}
We can see that $f_{ij}' \in P_k$ for all $k \geq 1$.   Note that $f_{ij}'$ is nonzero on the Hamiltonian cycle $x$ if
and only if in $x$, $i$ is adjacent to two vertices, neither of which is $j$; i.e. if and only if $i$ is not adjacent to $j$ in $x$.
  Also note that if $f_{ij}'(x) \not=0$, there is exactly one $\Gamma_x \in X_{ij}$ such that $g_{\Gamma_x}(x)
\not=0$.   Thus, we can see that $f_{ij}'$ is a constant multiple of the linear function $1-x_{ij}$, which corresponds to
a facet defined by the right hand side of equation \eqref{facet1}. 

Hence, we have shown that the facets defined by the left and right hand sides of equation \eqref{facet1} are on the
boundary of $P_k$ for $k \geq 1$.  

Suppose we have some $U \subset V$ with $\emptyset \not= U \not= V$, $|U| \leq 2k$. For each $i < |U|$, let 
\begin{align*}
X_i = &\{\Gamma : \Gamma \text{ corresponds to a path of length } i+1  \text{ with endpoints not in } U \\
& \text{ and } i \text{ vertices in } U\}
\end{align*}
Recall that $g_\Gamma(x)$, $\Gamma \in X_i$ takes on two values; 0 if $x$ does not contain $\Gamma$ and a positive number depending only on the size and type of partition corresponding to $\Gamma$ if $x$ does contain $\Gamma$.  For each $i<|U|$ let $c_i$ be a constant such that  if $\Gamma_i \in X_i$ and $x_i$ is a  Hamiltonian cycle such that $g_{\Gamma_i}(x_i) \not= 0$, then 
\begin{equation*}
c_ig_{\Gamma_i}(x_i)=2
\end{equation*}
and consider the function
\begin{equation*}
h_U = \sum_{i=1}^{|U|-1}\frac{|U|-i}{|U|}c_i\sum_{\Gamma \in X_i}g_{\Gamma}
\end{equation*}
Let $x$ be any Hamiltonian cycle.  Note that there will be an even number of edges, say $2\ell$ edges, from $U$ to $V-U$ in
$x$.   These will correspond to $\ell$ paths with all vertices except the endpoints in $U$.  If we sum over  those $\ell$ paths
the  number of vertices that each of the paths has in $U$, we will get $|U|$.  Those $\ell$ paths will correspond to the
only  $\Gamma \in X_i$ such that $g_{\Gamma}(x) \not= 0$.  Thus, we can see that 
\begin{equation*}
h_U(x) = \sum_{\substack{j \in U \\ i \in V-U}} x_{ij} -2
\end{equation*}
Hence, scaling $h_U$ so that we have a convex combination of the $g_\Gamma$s, we can see that for $|U| \leq 2k$, the facets corresponding to \eqref{facet2} are on the boundary of $P_k$.

\section{Proofs of Lemmas}\label{proofs}

Before we prove the Lemmas from section \ref{ProjectionConstruction}, we need one more Lemma.

\begin{Lemma}\label{pathcount} Let $m_1, \dotsc, m_p$ be nonnegative integers.  Then the number of ways of picking $m_1$ paths of length 1, $m_2$ paths of length 2, $\dotsc, m_p$ paths of length $p$ all from a path of $n$ vertices such that none of the chosen paths intersect is 
\begin{equation*}
\frac{\left( n-\sum_{i=1}^p im_i\right)_{\left(\sum_{i=1}^p m_i\right)}}{\prod_{i=1}^p m_i!}
\end{equation*}
\end{Lemma}

\begin{proof}[Proof of Lemma \ref{pathcount}]
The value above gives the number of ways of coloring $\sum_{i=1}^p m_i$ of the numbers from $1$ to $n-\sum_{i=1}^pim_i$ so that $m_i$ numbers are colored with color $i$.  We will construct a unique set of paths as required by the Lemma from each such coloring, and show that any set of paths can be obtained by a coloring.  We will present this bijection explicitly in the cases where we have exactly 1 or 2 total paths, and the cases where there are more paths will follow inductively.  

Suppose we have a path $P$ with $n$ vertices in it; number the vertices $1, 2, \dotsc, n$ so that 1 and $n$ are endpoint vertices, and $i$ is adjacent to $i-1$ and $i+1$ for $1 < i < n$.  For ease, we will think of $P$ as if we can visualize it horizontally, so that vertex 1 is to the left of vertex 2, which is to the left of vertex 3, etc.

Fix a number $j, 1 \leq j \leq n-1$.  Consider a set containing the numbers from 1 to $n-j$ such that one of these numbers is colored.  To this coloring we associate the path of length $j$ in $P$ which has its leftmost vertex located at the colored number.  In other words, if $i$ is the number which is colored in our set, we associate this coloring to the path of length $j$ in $P$ whose leftmost vertex is the vertex labeled $i$.  Since there are $j$ vertices to the right of the leftmost vertex in our path of length $j$, we see that this is a 1-1 correspondence between paths of length $j$ within a path of $n$ vertices and a coloring of one number from the set $\{1, 2, \dots, n-j\}$.

Now we consider a set containing the numbers from 1 to $n-j-i$ such that one of these numbers, say $n_i$, is colored with color $i$, and one, say $n_j$, is colored with color $j$.  WLOG, suppose that $n_i<n_j$.  If $n_j > n_i+i$, then to this coloring we associate the paths in $P$ of length $i$ and $j$ such that the path of length $i$ has the vertex labeled $n_i$ as its leftmost vertex and the path of length $j$ has the vertex labeled $n_j$ as its leftmost vertex.  A nice way to think of this is: given a coloring, we look at the path in $P$ of length $i$ with $n_i$ as its leftmost vertex and the path in $P$ of length $j$ with $n_j$ as its leftmost vertex.  If these paths do not intersect (i.e. if $n_j > n_i+i$), then those are the paths we associate to the coloring.  Otherwise, if they do intersect, we must move one of them.

Specifically, if $n_j = n_i+\ell$ where $\ell \leq i$, the paths intersect.  In this case, we move the path of length $j$ so that its leftmost vertex is $n-(j+\ell-1)$.  In other words, if there are $\ell'$ vertices between $n_i$ and $n_j$, we move the path of length $j$ in $P$ so that if $n_j'$ is the rightmost vertex of the moved path, there are exactly $\ell'$ vertices in $P$ to the right of $n_j'$.

This association of colorings of two numbers between $1$ and $n-j-i$ to placing paths of length $i$ and $j$ in $P$ is clearly reversible.  Specifically, suppose we have a path of length $i$ in $P$ with leftmost vertex $n_i$ and a path of length $j$ in $P$ with leftmost vertex $n_j$ which do not intersect.  If $n_i$ and $n_j$ are both less than $n-i-j$, the associated coloring is: color $n_i$ color $i$ and color $n_j$ color $j$.  Otherwise, one of them is larger than $n-i-j$ (note that not both of them can be larger, because the two paths do not intersect).  WLOG, say $n_j > n-i-j$.  Then there are $\ell <i$ vertices to the right of the path of length $j$.  So we color vertex $n_i$ with color $i$, and vertex $n_i+ \ell + 1$ with color $j$.  

The association between colorings and path placings when we have more than 2 paths to place is a simple iteration of the procedure when there are 1 or 2 paths to place.  Specifically, suppose we have a coloring of the numbers between 1 and $n-\sum_{i=1}^p im_i$ so that $m_i$ numbers are colored with color $i$.  From this coloring, we get a path placement as follows: let $n_i$ be the smallest colored number, colored with color $i$.  Place a path of length $i$ in $P$ with its leftmost vertex at vertex $n_i$.  Remove $n_i$ from the list of colored numbers.  Label $n$ as the ``last'' vertex in $P$.  Now at each step remaining, we do the following: find $n_j$ the smallest of the remaining colored numbers.  If vertex $n_j$ does not intersect one of the previously placed paths, place a path of length $j$ with leftmost vertex at $n_j$.  Then remove $n_j$ from the list of colored numbers and label $n-i$ as the ``last'' vertex (where $i$ is such that the closest colored vertex to the left of $n_j$ is colored with color $i$).  Otherwise, say $n_\ell$ is the leftmost vertex of a previously placed path of length $\ell$ and $n_j \leq n_\ell + \ell$.  Let $n_{last}$ be the vertex labeled ``last''  and  $n_j' = n_{last}-j+(n_\ell-n_j)-1$.  Then we place a path of length $j$ with leftmost vertex $n_j'$.  Remove $n_j$ from the list of colored numbers, and label $n_j'-1$ as the ``last'' vertex.  Continue inductively.  From the way that the ``last'' vertex is changed and from the fact that this is a 1-1 correspondence when we have 1 or 2 paths, we can see that this general construction works inductively.  Hence, we have proven the Lemma.

\end{proof}

\begin{proof}[Proof of Lemma \ref{y_0}]  Let $\pi$ be a partition of $k$ with $m_i$ parts of size $i$, so that $\sum_{i} i m_i = k$.  Let $m = \sum_{i} m_i$ be the number of parts in $\pi$.  Let $y$ be  a hamiltonian cycle which contains the edge $\{s,t\}$.  Recall
equation \eqref{explicitT}
\begin{align}
T_\pi(g_{st})(y)  = &\frac{|X|(n-1)}{a_\pi 2 |B_\pi| a_\pi} \biggl(|\{\Gamma \subset y : \Gamma \text{ has partition } \pi, \text{ and } \{s,t\} \in \Gamma\}| \notag \\
& \cdot 2^{m-1}(n-k-1)! \notag \\
+& |\{\Gamma \subset y : \Gamma \text{ has partition } \pi, \text{ and } s,t \text{ are each endpoints}  \notag \\
& \text{ of different paths in } \Gamma \}| \cdot 2^{m-2}(n-k-2)! \notag \\
+& |\{\Gamma \subset y : \Gamma \text{ has partition } \pi, \text{ and exactly one of } s,t  \text{ is an endpoint  } \notag \\ 
& \text{ of a path in } \Gamma \}| \cdot 2^{m-1}(n-k-2)! \notag \\
+& |\{ \Gamma \subset y : \Gamma \text{ has partition } \pi \text{ and } s,t  \text{ are not in } \Gamma \}|  \notag \\
& \cdot 2^m(n-k-2)! \biggr) \notag
\end{align}
By the same argument as in the proof of Theorem \ref{project1}, we can see that in this case 
\begin{equation*}
|B_\pi| = \frac{n_{(\sum_i (i+1)m_i)}}{2^{m} \prod_i m_i!} = \frac{n_{(k+m)}}{2^{m} \prod_i m_i!}
\end{equation*}
Also, using Lemma \ref{paths}, we can calculate 
\begin{equation*}
a_\pi = 2^{m -1}(n-k-1)!
\end{equation*}
Thus, we can see that 
\begin{align}
T_\pi(g_{st})(y)  = &\frac{(n-1)_{(k)}(n-1) \prod_i m_i!}{2n_{(k+m)}} \frac{1}{2^{m-1}(n-k-1)!}\biggl(|\{\Gamma \subset y : \Gamma \text{
has } \notag \\
& \text{partition } \pi,  \text{ and } \{s,t\} \in \Gamma\}|  \cdot 2^{m-1}(n-k-1)! \notag \\
+& |\{\Gamma \subset y : \Gamma \text{ has partition } \pi, \text{ and } s,t \text{ are each endpoints}  \notag \\
& \text{ of different paths in } \Gamma \}| \cdot 2^{m-2}(n-k-2)! \notag \\
+& |\{\Gamma \subset y : \Gamma \text{ has partition } \pi, \text{ and exactly one of } s,t  \text{ is an endpoint  } \notag \\ 
& \text{ of a path in } \Gamma \}| \cdot 2^{m-1}(n-k-2)! \notag \\
+& |\{ \Gamma \subset y : \Gamma \text{ has partition } \pi \text{ and } s,t  \text{ are not in } \Gamma \}|  \notag \\
& \cdot 2^m(n-k-2)! \biggr) \notag \\
 \geq & \frac{(n-1)_{(k)}(n-1) \prod_i m_i!}{2n_{(k+m)}} \frac{1}{2^{m-1}(n-k-1)!}\biggl(|\{\Gamma \subset y : \Gamma \text{ has
} \notag \\
&\text{partition } \pi, \text{ and } \{s,t\} \in \Gamma\}|  \cdot 2^{m-1}(n-k-1)! \notag \\
 +& |\{ \Gamma \subset y : \Gamma \text{ has partition } \pi \text{ and } s,t  \text{ are not in } \Gamma \}|  \notag \\
 & \cdot 2^m(n-k-2)! \biggr) \notag 
\end{align}
We can count $|\{\Gamma \subset y : \Gamma \text{ has partition } \pi, \text{ and } \{s,t\} \in \Gamma\}|$ as follows: suppose that $m_i \not=0$ (this is true for at least some $i$).  Mark a path of length $i$ arbitrarily in the Hamiltonian cycle $y$.  The number of ways that the remaining paths can be chosen, using Lemma \ref{pathcount}, is 
\begin{equation*}
\frac{(n-(i+1)-\sum_{j \not=i}jm_j-(m_i-1)i)_{((m_i-1)+\sum_{j \not=i}m_j)}}{(m_i-1)!\prod_{j\not=i} m_j!} = \frac{(n-k-1)_{(m-1)}}{(m_i-1)!\prod_{j\not=i}m_j!}
\end{equation*}
Then there are $k$ ways of rotating the cycle $y$ cyclically so that the edge $\{s,t\}$ lies in one of our chosen paths.  Since the $m_i$ paths of length $i$ are indistinguishable, we need to divide by $m_i$ in order to not overcount.  Thus, we have found
\begin{equation*}
|\{\Gamma \subset y : \Gamma \text{ has partition } \pi, \text{ and } \{s,t\} \in \Gamma\}| = \frac{k(n-k-1)_{(m-1)}}{\prod_{j}m_j!}
\end{equation*}
Using Lemma \ref{pathcount}, we can calculate
\begin{equation*}
 |\{ \Gamma \subset y : \Gamma \text{ has partition } \pi \text{ and } s,t  \text{ are not in } \Gamma \}|  = \frac{(n-2-k)_{(m)}}{\prod_j m_j!}
 \end{equation*}
 Thus, since we assume $k < n^{1/3}$, we have shown that 
 \begin{align}
T_\pi(g_{st})(y)  \geq & \frac{(n-1)_{(k)}(n-1) \prod_i m_i!}{2n_{(k+m)}} \frac{1}{2^{m-1}(n-k-1)!}  \notag \\ 
&\cdot\biggl( \frac{k(n-k-1)_{(m-1)}}{\prod_{j}m_j!}   \cdot 2^{m-1}(n-k-1)!\notag \\
 &+\frac{(n-2-k)_{(m)}}{\prod_j m_j!} \cdot 2^m(n-k-2)! \biggr) \notag \\
 =& \frac{k+2}{2} +  O\left(\frac{1}{n^{2/3}}\right) \notag
\end{align}

\end{proof}

\begin{proof}[Proof of Lemma \ref{same}]

Looking at equation \eqref{explicitT}, we see that we can think of $T_\pi(g_{st})(y)$ as a sum over all $\Gamma \subset y$
 with partition type $\pi$, each $\Gamma$ contributing a certain amount (maybe 0).    We need to show that $T_\pi(g_{st})(y)$ is
 the same for all Hamiltonian cycles $y$ such that the distance between $s$ and $t$ in $y$ is at least $k$.   We will do
this by showing that if $y$ is a Hamiltonian cycle with $d$ vertices between $s$ and $t$, $k \leq d < \frac{n}{2}$,  
then for a Hamiltonian cycle $y'$ with $d+1$ vertices between $s$ and $t$ we have $T_\pi(g_{st})(y) =T_\pi(g_{st})(y')$.
  Since the only thing that affects the value of $T_\pi(g_{st})(y)$ is the number of vertices between $s$ and $t$ in
$y$,  we can WLOG consider the following two cases 
\begin{align*}
y &= (1, 2, \dotsc, n) &s = 1, \quad t = d  \\
y' &= (1, 2,  \dotsc, n) &s'=1, \quad t'=d+1 
\end{align*}
where $d \geq k+2$.  (Hence we will show $T_\pi(g_{st})(y) = T_\pi(g_{s't'})(y')$.)

Consider any $\Gamma \subset y$ with partition type $\pi$.  Consider the exact same $\Gamma \subset y'$.   Note that,
since $\pi$ is a partition of $k$ and there are at least $k+1$ edges from vertex 1 to vertex $d$ (at least $k+2$ edges
from vertex 1 to vertex $d+1$), $\Gamma$  cannot have
a path in it connecting vertex $1$ to vertex $d$ (connecting vertex $1$ to vertex $d+1$).  Suppose that one of the following is true: 
\begin{enumerate}
\item $d$ and $d+1$ are each endpoints of some path in $\Gamma$
\item Neither $d$ nor $d+1$ appears in $\Gamma$
\item Both $d$ and $d+1$ are in the middle (\textit{not} an endpoint) of a path in $\Gamma$
\end{enumerate}
Then, in looking at equation \eqref{explicitT}, we see that this particular $\Gamma$ contributes the same amount in $T_\pi(g_{st})(y)$ as in $T_\pi(g_{s't'})(y')$.  The only cases where $\Gamma$ contributes differing amounts in $T_\pi(g_{st})(y)$ and $T_\pi(g_{s't'})(y')$ are:
\begin{enumerate}
\item One of $d$ or $d+1$ is an endpoint of a path in $\Gamma$, and the other does not appear in $\Gamma$
\item One of $d$ or $d+1$ is an endpoint of a path in $\Gamma$, and the other is in the middle (\textit{not} and endpoint) of a path in $\Gamma$.
\end{enumerate}

If $\Gamma$ lies in one of those two cases, it contributes a different amount in $T_\pi(g_{st})(y)$ versus in $T_\pi(g_{s't'})(y')$. 
Thus, to each $\Gamma$ in one of those
 two cases, we will associate a unique $\Gamma'$
(also in one of those two cases) for which the contribution of $\Gamma$ in $T_\pi(g_{st})(y)$ is equal to the
contribution of $\Gamma'$ in
$T_\pi(g_{s't'})(y')$, and the contribution of $\Gamma'$ in $T_\pi(g_{st})(y)$ is equal to the contribution of $\Gamma$ in $T_\pi(g_{s't'})(y')$. 
We will show that if the partner to $\Gamma$ under this association is $\Gamma'$, then the partner to $\Gamma'$ under 
this association is $\Gamma$.  Once we have this, we will be done.

Given some $\Gamma$ in one of the two cases above, consider vertices $d-1$ and $d+2$.  If vertex $d-1$ is not connected
to vertex $d-2$ \textit{and} vertex $d+2$ is not connected to vertex
$d+3$,  we map $\Gamma \mapsto \Gamma'$, where $\Gamma'$ leaves all paths in $\Gamma$ untouched, except for the
paths from $d-1$ to $d+2$ which it  ``reflects'' 
about the line between $d$ and $d+1$, as demonstrated in the following picture:
\begin{gather*}
\circ_{d-1} - \circ_d \quad  \circ_{d+1} \quad  \circ_{d+2} \quad \quad \leftrightarrow \quad \quad \circ_{d-1} \quad  \circ_d \quad \circ_{d+1} - \circ_{d+2} \\
\text{and} \\
\circ_{d-1} - \circ_d -  \circ_{d+1} \quad  \circ_{d+2} \quad \quad \leftrightarrow \quad \quad \circ_{d-1} \quad  \circ_d - \circ_{d+1} - \circ_{d+2}
\end{gather*}
Since vertex $d-1$ is not connected to vertex $d-2$ and vertex $d+2$ is not connected to vertex $d+3$, this action preserves
the partition type so that we obtain a partner $\Gamma'$ again of partition type $\pi$.  It is also clear that this action
indeed  produces a $\Gamma'$ such that the contribution of $\Gamma$ in $T_\pi(g_{st})(y)$ is equal to the contribution
 of $\Gamma'$ in $T_\pi(g_{s't'})(y')$, and the contribution of $\Gamma'$ in $T_\pi(g_{st})(y)$ is equal to the
contribution of $\Gamma$ in $T_\pi(g_{s't'})(y')$.
 We also note that $\Gamma'$ is mapped to $\Gamma$ under this action.

If $d-1$ was connected to $d-2$ or if $d+2$ was connected to $d+3$, then check to see if $d-2$ is connected to $d-3$ and if $d+3$ is connected to $d+4$.
  If $d-2$ is not connected to $d-3$ \textit{and} $d+3$ is not connected to $d+4$, we can do the same ``reflecting'' action, this time between the paths from $d-2$ to $d+3$.  If not, continue looking for the first place where we can reflect.

Firstly we note that the first place to reflect is well-defined, and that if this action associates $\Gamma$ to
$\Gamma'$, it associates $\Gamma'$ to $\Gamma$. 
Secondly, we note that we will find a
``first place to reflect'' before getting down to vertex 1.   This is because $d \geq k+2$ and there are only $k$ edges
in $\Gamma$.  Finally, it is clear that this action indeed produces a  $\Gamma'$ such that the contribution of $\Gamma$
in $T_\pi(g_{st})(y)$ is equal to the contribution of $\Gamma'$ in $T_\pi(g_{s't'})(y')$, and the contribution of $\Gamma'$ in
 $T_\pi(g_{st})(y)$ is equal to the contribution of $\Gamma$ in $T_\pi(g_{s't'})(y')$.  Thus, we have proven the Lemma.

\end{proof}

The ideas of ``reflecting'' in the proof of this lemma will be also come into play in the proofs of Lemmas
\ref{change_sequence}-\ref{ksnake}.

\begin{proof}[Proof of Lemma \ref{change_sequence}]

The proof of Lemma \ref{same} shows us exactly how to prove Lemma \ref{change_sequence}.  Namely, the reason that
$T_\pi(g_{st})$ takes on different values for Hamiltonian cycles having $i$ and $i+1$ vertices between $s$ and $t$ ($i<k$)
is  because the association from $\Gamma$ to $\Gamma'$ may not work.  In particular, let $\pi$ be the partition $(\underbrace{1, 1, \dotsc, 1}_{k \; \text{ 1s}})$ and
\begin{align*}
y &= (1, 2, \dotsc, n) &s = 1, \quad t = i+2  \\
y' &= (1, 2,  \dotsc, n) &s'=1, \quad t'=i+3 
\end{align*}
where $1 \leq i <k-1$.  (Then $y$  has $i$ vertices between $s$ and $t$, $y'$ has $i+1$ vertices between $s$ and $t$).  Then each path subset $\Gamma$ of partition type $\pi$ can be associated to a $\Gamma'$ just as in the proof of Lemma \ref{same}, unless we have something like the following:
\begin{gather}
\circ_{1} - \circ_2 \quad  \circ_{3} -  \circ_{4} \quad   \dotso \quad  \circ_{i+1} - \circ_{i+2}  \quad \circ_{i+3} \quad \circ_{i+4} - \circ_{i+5} \dotso \circ_{2i+2} - \circ_{2i+3} \label{first1}\\
\text{or} \notag \\
\circ_{1} \quad \circ_2 -  \circ_{3} \quad  \circ_{4} - \circ_{5} \dotso  \circ_{i}-\circ_{i+1} \quad \circ_{i+2}  \circ_{i+3} - \circ_{i+4}\dotso \circ_{2i+3} - \circ_{2i+4} \label{first2}
\end{gather}
in the case of $i$ even or 
\begin{gather}
\circ_{1} - \circ_2 \quad  \circ_{3} -  \circ_{4} \quad   \dotso \quad  \circ_{i} - \circ_{i+1}  \quad \circ_{i+2} \quad \circ_{i+3} - \circ_{i+4} \dotso \circ_{2i+2} - \circ_{2i+3} \label{first3}\\
\text{or} \notag \\
\circ_{1} \quad \circ_2 -  \circ_{3} \quad  \circ_{4} - \circ_{5} \dotso  \circ_{i+1}-\circ_{i+2} \quad \circ_{i+3} \quad \circ_{i+4}-\circ_{i+5} \dotso \circ_{2i+3} - \circ_{2i+4} \label{first4}
\end{gather}
in the case of $i$ odd.

Thus, we can see that the \textit{difference} between $T_\pi(g_{st})(y)$ and $T_\pi(g_{s't'})(y')$ is  simply the
difference in the contribution of each path subset $\Gamma$ which does not have a valid partner to which it can map (i.e., if $\Gamma$ corresponds to one of the above cases).  
 Thus, in looking at equation \eqref{explicitT} and using the notation of Lemma \ref{change_sequence}, we can see that
$T_\pi(g_{st})(y)-T_\pi(g_{s't'})(y')$ is

\begin{align*}
\frac{|X|(n-1)}{a_\pi 2 |B_\pi| a_\pi} \biggl(
 & |\{\Gamma \text{ in case of \eqref{first2} and } \Gamma \text{ contains edge} \{1,n\} \}| \\
&\cdot 2^{k-2}(n-k-2)!  \\
+ & |\{\Gamma \text{ in case of \eqref{first2} and } \Gamma \text{ does not contain edge } \{1,n\}  \\
& \text{or } \Gamma \text{ in case of \eqref{first1}}\}| \cdot 2^{k-1}(n-k-2)!  \\
-& |\{\Gamma \text{ in case of \eqref{first1}} \}| \cdot 2^{k-2}(n-k-2)!  \\
- & |\{\Gamma \text{ in case of \eqref{first2} and } \Gamma \text{ contains edge } \{1,n\} \}| \cdot 2^{k-1}(n-k-2)!  \\
-&|\{\Gamma \text{ in case of \eqref{first2} and } \Gamma \text{ does not contain edge } \{1,n\} \}| \\
 & \cdot 2^k(n-k-2)! \biggr)
\end{align*}
when $i$ is even and 

\begin{align*}
\frac{|X|(n-1)}{a_\pi 2 |B_\pi| a_\pi} \biggl(
 & |\{\Gamma \text{ in case of \eqref{first3}} \}| \cdot 2^{k-2}(n-k-2)!  \\
+ & |\{\Gamma \text{ in case of \eqref{first4} and } \Gamma \text{ contains edge } \{1,n\} \}| \cdot 2^{k-1}(n-k-2)!  \\
+ & |\{\Gamma \text{ in case of \eqref{first4} and } \Gamma \text{ does not contain edge } \{1, n\} \}| \\
& \cdot  2^k(n-k-2)! \\
-& |\{\Gamma \text{ in case of \eqref{first4} and } \Gamma \text{ contains edge } \{1, n\}\}|  \cdot 2^{k-2}(n-k-2)!  \\
- & |\{\Gamma \text{ in case of \eqref{first4} and } \Gamma \text{ does not contain edge } \{1,n\}  \\
& \text{or } \Gamma \text{ in case of \eqref{first3}} \}|\cdot 2^{k-1}(n-k-2)!  \\\biggr)
\end{align*}
when $i$ is odd.

Using Lemma \ref{pathcount}, we can actually calculate these differences.  In the case of $i$ even, we have
\begin{align}
\frac{|X|(n-1)}{a_\pi 2 |B_\pi| a_\pi} \biggl(
& \frac{(n-(2i+5)-(k-i-2))_{(k-i-2)}}{(k-i-2)!}2^{k-2}(n-k-2)! \notag \\
+&\left(\frac{(n-(2i+4)-(k-i-1))_{(k-i-1)}}{(k-i-1)!} \right. \notag \\
+& \left. \frac{(n-(2i+3) -(k-i-1))_{(k-i-1)}}{(k-i-1)!}\right) \cdot 2^{k-1}(n-k-2)! \notag \\
-&\frac{(n-(2i+3) -(k-i-1))_{(k-i-1)}}{(k-i-1)!}2^{k-2}(n-k-2)!  \notag \\
-& \frac{(n-(2i+5)-(k-i-2))_{(k-i-2)}}{(k-i-2)!}2^{k-1}(n-k-2)! \notag \\
-& \frac{(n-(2i+4)-(k-i-1))_{(k-i-1)}}{(k-i-1)!}2^k(n-k-2)!  \biggr) \label{bigmess1}
\end{align}
and in the case of $i$ odd we have
\begin{align}
\frac{|X|(n-1)}{a_\pi 2 |B_\pi| a_\pi} \biggl(
& \frac{(n-(2i+3)-(k-i-1))_{(k-i-1)}}{(k-i-1)!}2^{k-2}(n-k-2)! \notag \\
+& \frac{(n-(2i+5)-(k-i-2))_{(k-i-2)}}{(k-i-2)!}2^{k-1}(n-k-2)! \notag \\
+& \frac{(n-(2i+4)-(k-i-1))_{(k-i-1)}}{(k-i-1)!}2^k(n-k-2)! \biggr) \notag \\
-& \frac{(n-(2i+5)-(k-i-2))_{(k-i-2)}}{(k-i-2)!}2^{k-2}(n-k-2)! \notag \\
-& \left(\frac{(n-(2i+4)-(k-i-1))_{(k-i-1)}}{(k-i-1)!}\right. \notag \\
-&\left. \frac{(n-(2i+3)-(k-i-1))_{(k-i-1)}}{(k-i-1)!}\right) \cdot 2^{k-1}(n-k-2)! \biggr) \label{bigmess2}
\end{align}

Thus, to calculate these values, all we have left is to compute the values of $|X|$, $a_\pi$, and $|B_\pi|$.  We have already identified that $|X|$, the number of Hamiltonian cycles in the complete graph on $n$ vertices $K_n$, is $\frac{(n-1)!}{2}$.  Recall that $|B_\pi|$ is the number of path subsets of partition type $\pi$.  By the same argument used in the proof of Theorem \ref{project1}, we have:
\begin{equation*}
|B_\pi| = \frac{n_{(2k)}}{2^kk!}
\end{equation*}

Finally, recall that $a_\pi$ is the number of Hamiltonian cycles containing all edges in a path subset $\Gamma$ of partition type $\pi$.  Then from Lemma \ref{paths}, we know that for $\pi=(\underbrace{1,1, \dotsc, 1}_{k \text{ 1s}})$, $a_\pi = 2^{k-1}{(n-k-1)!}$.  Plugging all of these into equations \eqref{bigmess1} and \eqref{bigmess2}, we find that for $i$ even we have

\begin{align*}
\frac{(n-1)_{(k)}(n-1)k!}{ n_{(2k)}(n-k-1)4} \biggl(&\frac{(n-(2i+5)-(k-i-2))_{(k-i-2)}}{(k-i-2)!}  \\
&+\frac{(n-(2i+4)-(k-i-1))_{(k-i-1)}}{(k-i-1)!}2  \\
&+  \frac{(n-(2i+3) -(k-i-1))_{(k-i-1)}}{(k-i-1)!}2  \\
&-\frac{(n-(2i+3) -(k-i-1))_{(k-i-1)}}{(k-i-1)!} \\
&- \frac{(n-(2i+5)-(k-i-2))_{(k-i-2)}}{(k-i-2)!}2  \\
&- \frac{(n-(2i+4)-(k-i-1))_{(k-i-1)}}{(k-i-1)!}4 \biggr)  \\
=\frac{(n-1)_{(k)}(n-1)k!}{ n_{(2k)}(n-k-1)4}  \biggl(&-2\frac{(n-k-i-3)_{(k-i-1)}}{(k-i-1)!}   \\
&+ \frac{(n-k-i-2)_{(k-i-1)}}{(k-i-1)!}  \\
&-\frac{(n-k-i-3)_{(k-i-2)}}{(k-i-2)!} \biggr)  
\end{align*}
\begin{equation}
= -\frac{(n-1)(n-2k)(n-2k-1)k_{(i+1)}}{4(n-k-1)n(n-k-1)_{(i+2)}}\label{finalmess1}
\end{equation}

and for $i$ odd we have
\begin{align*}
\frac{(n-1)_{(k)}(n-1)k!}{n_{(2k)}(n-k-1)4} \biggl(
& \frac{(n-(2i+3)-(k-i-1))_{(k-i-1)}}{(k-i-1)!}  \\
&+ \frac{(n-(2i+5)-(k-i-2))_{(k-i-2)}}{(k-i-2)!}2  \\
&+ \frac{(n-(2i+4)-(k-i-1))_{(k-i-1)}}{(k-i-1)!}4  \\
&- \frac{(n-(2i+5)-(k-i-2))_{(k-i-2)}}{(k-i-2)!}  \\
&- \frac{(n-(2i+4)-(k-i-1))_{(k-i-1)}}{(k-i-1)!}2  \\
&- \frac{(n-(2i+3)-(k-i-1))_{(k-i-1)}}{(k-i-1)!} 2 \biggr)  \\
=\frac{(n-1)_{(k)}(n-1)k!}{n_{(2k)}(n-k-1)4} \biggl(
&  2\frac{(n-k-i-3)_{(k-i-1)}}{(k-i-1)!}  \\
& - \frac{(n-k-i-2)_{(k-i-1)}}{(k-i-1)!}  \\
& + \frac{(n-k-i-3)_{(k-i-2)}}{(k-i-2)!} \biggr)   
\end{align*}
\begin{equation}
= \frac{(n-1)(n-2k)(n-2k-1)k_{(i+1)}}{4(n-k-1)n(n-k-1)_{(i+2)}}\label{finalmess2}
\end{equation}
and we have finished our proof.

\end{proof}

\begin{proof}[Proof of Lemma \ref{k,1snakes}]  
Recall that here $\pi = (k-1,1)$ and $k \geq 3$.  In this case, we can actually calculate the value of $T_\pi(g_{st})$ on any Hamiltonian cycle.
 Recall that
\begin{align}\label{k,1snakes_eqn}
T_\pi(g_{st})(y)  = &\frac{|X|(n-1)}{a_\pi 2 |B_\pi| a_\pi} \biggl(|\{\Gamma \subset y : \Gamma \text{ has partition } \pi, \text{ and } \{s,t\} \in \Gamma\}| \notag \\
& \cdot 2^{m-1}(n-k-1)! \notag \\
+& |\{\Gamma \subset y : \Gamma \text{ has partition } \pi, \text{ and } s,t \text{ are each endpoints}  \notag \\
& \text{ of different paths in } \Gamma \}| \cdot 2^{m-2}(n-k-2)! \notag \\
+& |\{\Gamma \subset y : \Gamma \text{ has partition } \pi, \text{ and exactly one of } s,t  \text{ is an endpoint  } \notag \\ 
& \text{ of a path in } \Gamma \}| \cdot 2^{m-1}(n-k-2)! \notag \\
+& |\{ \Gamma \subset y : \Gamma \text{ has partition } \pi \text{ and } s,t  \text{ are not in } \Gamma \}|  \notag \\
& \cdot 2^m(n-k-2)! \biggr)
\end{align}
Using the same argument used in the proof of Theorem \ref{project1}, we can calculate 
\begin{equation*}
|B_\pi| = \frac{n_{(k+2)}}{4}
\end{equation*}
From Lemma \ref{paths} we can calculate
\begin{equation*}
a_\pi = 2(n-k-1)!
\end{equation*}
 Suppose that $y_i$ has $i$ vertices between $s$ and $t$ for $1 \leq i \leq k-2$.  Then, using Lemma \ref{pathcount}
and counting all the ways that the different intersection patterns described in \eqref{k,1snakes_eqn} can occur, we calculate
\begin{align*}
T_\pi(g_{st})(y_i) =& \frac{4(n-1)(n-1)!}{16n_{(k+2)}(n-k-1)!(n-k-1)!}(4(n-k-2)!+(2(n-k-3) \\
&+2(n-k-i-2)+2(n-k-i-1))2(n-k-2)! \\
&+((i-1)(n-k-i-1)+(n-k-i-2)(n-k-i-3))\\
& \cdot 4(n-k-2)!) \\
=& \frac{n-1}{2n(n-k-1)}\frac{2n^2-(4k+2i+6)n+2k^2+2ik+6i+6k+4}{n-k-1}
\end{align*}
Suppose $y_{k-1}$ has $k-1$ vertices between $s$ and $t$.  Then we calculate 
\begin{align*}
T_\pi(g_{st})(y_{k-1}) =& \frac{4(n-1)(n-1)!}{16n_{(k+2)}(n-k-1)!(n-k-1)!}(6(n-k-2)!+(2(n-k-3) \\
&+2(n-k-2) +2(n-2k)+2(n-2k-1))2(n-k-2)! \\
&+((k-2)(n-2k)+(n-2k-1)(n-2k-2))4(n-k-2)!) \\
=& \frac{n-1}{2n(n-k-1)}\frac{2n^2-(6k+2)n+4k^2+8k-5}{n-k-1}
\end{align*}
Finally, suppose that $y_k$ has $k$ vertices between $s$ and $t$.  Then we calculate 
\begin{align*}
T_\pi(g_{st})(y_{k}) =& \frac{4(n-1)(n-1)!}{16n_{(k+2)}(n-k-1)!(n-k-1)!}(8(n-k-2)!+(4(n-k-3) \\
&+4(n-2k-1))2(n-k-2)! +((k-1)(n-2k-1) \\
&+(n-2k-2)(n-2k-3)+(n-k-3))4(n-k-2)!) \\
=& \frac{n-1}{2n(n-k-1)}\frac{2n^2-(6k+2)n+4k^2+8k-4}{n-k-1}
\end{align*}
The proof now follows.
\end{proof}

\begin{proof}[Proof of Lemma \ref{ksnake}]  
Recall that here $\pi = (k)$.  Again we can calculate the value of $T_\pi(g_{st})$ on any Hamiltonian cycle.  Suppose that $y_i$ has $i$ vertices between $s$ and $t$ for $1 \leq i \leq k-1$.  Recall that
\begin{align}\label{ksnake_eqn}
T_\pi(g_{st})(y)  = &\frac{|X|(n-1)}{a_\pi 2 |B_\pi| a_\pi} \biggl(|\{\Gamma \subset y : \Gamma \text{ has partition } \pi, \text{ and } \{s,t\} \in \Gamma\}| \notag \\
& \cdot 2^{m-1}(n-k-1)! \notag \\
+& |\{\Gamma \subset y : \Gamma \text{ has partition } \pi, \text{ and } s,t \text{ are each endpoints}  \notag \\
& \text{ of different paths in } \Gamma \}| \cdot 2^{m-2}(n-k-2)! \notag \\
+& |\{\Gamma \subset y : \Gamma \text{ has partition } \pi, \text{ and exactly one of } s,t  \text{ is an endpoint  } \notag \\ 
& \text{ of a path in } \Gamma \}| \cdot 2^{m-1}(n-k-2)! \notag \\
+& |\{ \Gamma \subset y : \Gamma \text{ has partition } \pi \text{ and } s,t  \text{ are not in } \Gamma \}|  \notag \\
& \cdot 2^m(n-k-2)! \biggr)
\end{align}
Using the same argument used in the proof of Theorem \ref{project1}, we can calculate 
\begin{equation*}
|B_\pi| = \frac{n_{(k+1)}}{2}
\end{equation*}
From Lemma \ref{paths} we can calculate
\begin{equation*}
a_\pi = (n-k-1)!
\end{equation*}
Thus, using Lemma \ref{pathcount} and counting all the ways that the different intersection patterns described in \eqref{ksnake_eqn} can occur, we calculate
\begin{align*}
T_\pi(g_{st})(y_i) =& \frac{2(n-1)(n-1)!}{4n_{(k+1)}(n-k-1)!(n-k-1)!}(2(n-k-2)! \\
&+(n-k-i-2)2(n-k-2)!) \\
=& \frac{(n-1)(n-k-i-1)}{n(n-k-1)}
\end{align*}
Suppose $y_{k}$ has $k$ vertices between $s$ and $t$.  Then we calculate 
\begin{align*}
T_\pi(g_{st})(y_{k}) =&\frac{2(n-1)(n-1)!}{4n_{(k+1)}(n-k-1)!(n-k-1)!}(4(n-k-2)! \\
&+(n-2k-2)2(n-k-2)!) \\
=& \frac{(n-1)(n-2k)}{n(n-k-1)}
\end{align*}
The proof now follows.
\end{proof}

\section{Remaining Comments}

Let $\pi^*, \pi_\ell,$ and $\pi_{\ell'}$ be as in the proof of Theorem \ref{project1}.  We note that we proved Lemmas
\ref{k,1snakes} and \ref{ksnake} by calculating the values of $T_{\pi_\ell}(g_{st})$ and $T_{\pi_{\ell'}}(g_{st})$.  
One can use Lemma \ref{pathcount} to calculate the value of $T_{\pi^*}(g_{st})$ on Hamiltonian cycles having 1  vertex
between $s$ and $t$, and use Lemma \ref{change_sequence} to calculate the remaining values of $T_{\pi^*}(g_{st})$.  
Calculating values of $T_{\phi}(g_{st})$ for arbitrary partitions $\phi$ can be done using equation \eqref{explicitT} 
and Lemma \ref{pathcount}.  However, determining the leading terms of \eqref{explicitT} for arbitrary partitions and 
arbitrary Hamiltonian cycles is much more complicated when $k$ is not fixed.  It is a reasonable question whether one can
find  a convex combination of $T_\phi$ for other partitions $\phi$ which gives a better approximation, or which is valid for $k$ closer to $n$.

The role of the linear operators $T_\pi$ in the proof of Theorem \ref{project1} is to show that a  scaling of the set
$Q-\mathbbm{1}$ lies inside of the set $P_k \cap L -\mathbbm{1}$.  It is entirely possible that the  scaling factor we
achieve using these linear operators is not optimal; that we could scale $Q-\mathbbm{1}$ by a larger number  and have
it still lie inside $P_k\cap L -\mathbbm{1}$.

\bibliographystyle{amsplain}
\bibliography{TSPapprox}
\end{document}